\documentclass[twoside,10pt]{jmlr} % Anonymized submission

\usepackage{fourier}
\usepackage{eucal}

\usepackage{wrapfig}
\usepackage{caption}
\usepackage{subcaption}
\usepackage[utf8]{inputenc} % allow utf-8 input
\usepackage[T1]{fontenc}    % use 8-bit T1 fonts
\usepackage{hyperref}       % hyperlinks
\usepackage{url}            % simple URL typesetting
\usepackage{booktabs}       % professional-quality tables
\usepackage{amsfonts}       % blackboard math symbols
\usepackage{nicefrac}       % compact symbols for 1/2, etc.
\usepackage{microtype}      % microtypography
\usepackage[utf8]{inputenc}
\usepackage{amsmath}
\newcommand{\qed}{\hfill\ensuremath{\blacksquare}}%
\newtheorem{prop}{Proposition}
\newtheorem{assumption}{Assumption}
\usepackage{tcolorbox}
\usepackage{graphicx}
\usepackage{subcaption}

\newcommand{\X}{\mathcal{X}}
\renewcommand{\S}{\mathcal{S}}
\newcommand{\Z}{\mathcal{Z}}

\newcommand{\A}{\mathcal{A}}

\newcommand{\real}{\mathbb{R}}

\newcommand{\DD}[2]{D_\Phi\!\left(#1\middle\|#2\right)}
\newcommand{\OO}{\mathcal{O}}
\newcommand{\LL}{\mathcal{L}}

\newcommand{\II}[1]{\mathbb{I}_{\left\{#1\right\}}}
\newcommand{\PP}[1]{\mathbb{P}\left[#1\right]}

\newcommand{\EEpi}[1]{\mathbb{E}_\pi\left[#1\right]}

\def\argmax{\mathop{\mbox{ arg\,max}}}
\newcommand{\ra}{\rightarrow}

\newcommand{\taumix}{\tau_{\mbox{\scriptsize{mix}}}}
\newcommand{\iprod}[2]{\left\langle#1,#2\right\rangle}
\newcommand{\norm}[1]{\left\|#1\right\|}
\newcommand{\onenorm}[1]{\norm{#1}_1}
\newcommand{\twonorm}[1]{\norm{#1}_2}
\newcommand{\infnorm}[1]{\norm{#1}_\infty}

\newcommand{\ev}[1]{\left\{#1\right\}}
\newcommand{\pa}[1]{\left(#1\right)}
\newcommand{\bpa}[1]{\bigl(#1\bigr)}

\newcommand{\wh}{\widehat}
\newcommand{\wt}{\widetilde}

\newcommand{\transpose}{^\mathsf{\scriptscriptstyle T}}

\newcommand{\tD}{\wt{\Delta}}
\newcommand{\hz}{\wh{z}}

\newcommand{\bmu}{\overline{\mu}}
\newcommand{\bnu}{\overline{\stdist}}

\newcommand{\bu}{\overline{u}}
\newcommand{\by}{\overline{y}}

\usepackage{todonotes}
\definecolor{PalePurp}{rgb}{0.66,0.57,0.66}

\newcommand{\stdist}{d}
\newcommand{\hy}{\wh{y}}
\newcommand{\hu}{\wh{u}}
\newcommand{\hye}{\hy_\varepsilon}

\title[Faster saddle-point optimization for solving large-scale MDPs]{Faster saddle-point 
optimization for solving \\ large-scale Markov decision processes}

\author[Bas-Serrano and Neu]{\Name[{Joan~Bas-Serrano}]{Joan Bas-Serrano} 
\Email{joanbasserrano@gmail.com}\\
\addr Universitat Pompeu Fabra, Barcelona, Spain
\\
\Name[{Gergely~Neu}]{Gergely Neu} \Email{gergely.neu@gmail.com}\\
\addr Universitat Pompeu Fabra, Barcelona, Spain
 }

\begin{document}

\maketitle

\begin{abstract}
We consider the problem of computing optimal policies in average-reward Markov decision processes. 
This classical problem can be formulated as a linear program directly amenable to saddle-point 
optimization methods, albeit with a number of variables that is linear in the number of states. To 
address this issue, recent work has considered a linearly relaxed version of the resulting 
saddle-point problem. Our work aims at achieving a better understanding of this relaxed 
optimization problem by characterizing the conditions necessary for convergence to the 
optimal policy, and designing an optimization algorithm enjoying fast convergence rates that are
independent of the size of the state space. Notably, our characterization points out some potential 
issues with previous work.
\end{abstract}

\section{Introduction}
Computing optimal policies in Markov decision processes (MDPs) is one of the most important 
problems in sequential decision making and control \citep{Puterman1994}. Arguably, the most 
classical approach to solve this task is through the method of \emph{dynamic programming}, 
understood in this 
context as computing fixed points of certain operators \citep{Bel57,How60,Ber07:DPbookVol2}. The 
use 
and influence of dynamic-programming methods like value and policy iteration extend well beyond the 
world of decision and control theory, as the underlying ideas serve as foundations for most 
algorithms for \emph{learning} optimal policies in unknown MDPs: the setting of \emph{reinforcement 
learning} \citep{Sze10,SB18}. While being hugely successful, DP-based methods have the downside of 
being somewhat incompatible with classical machine-learning tools that are rooted in convex 
optimization. Indeed, most of the popular reductions of dynamic programming to (non-)convex 
optimization are based on heuristics that are not directly motivated by theory. Examples include 
the celebrated DQN approach of \citet{mnih2015human} that reduces value-function estimation to 
minimizing the ``squared Bellman error'', or the TRPO algorithm of \citet{SLAJM15} that reduces 
policy updates to minimizing a ``regularized surrogate objective''. While these methods can be 
justified to a certain extent, it is technically unknown if solving the resulting optimization 
problems actually leads to a desirable solution to the original sequential decision-making problem.

In this paper, we explore a family of methods that reduce MDP optimization to a form of convex 
optimization in a theoretically grounded way. Our starting point is an alternative approach based 
on linear programming (LP), first proposed roughly at the same time as the DP methods of 
\citet{Bel57,How60}: the idea of LP-based methods for sequential decision-making goes back to the 
works of \citet{Ghe60,Man60,Den70}. While LP-based methods seem to be more obscure in present 
day than DP methods, they have the clear advantage that they lead to an objective function directly 
amenable to modern large-scale optimization methods. Recent reinforcement-learning methods inspired 
by the LP perspective include policy-gradient and actor-critic methods \citep{SuMcSiMa99,KoTsit99} 
and various ``entropy-regularized'' learning algorithms (e.g., 
\citealp{peters10reps,ZiNe13,NJG17}). While 
these methods 
promise to directly tackle the policy-optimization problem through solving the underlying linear 
program, most of them still require the computation of certain value functions through dynamic 
programming.

In the present work, we argue for the viability of a method fully based on a form of convex 
optimization, rooted in the LP approach. Our approach is based on a \emph{bilinear saddle-point} 
formulation of the linear program, building on a well-known general equivalence between the two 
optimization problems. One particular advantage of this formulation is that it enables a 
straightforward form of dimensionality reduction of the original problem through a linear 
parametrization of the optimization variables, which provides a natural framework for studying 
effects of ``function approximation'' in the underlying policy optimization problem. Our main 
contribution regarding this setting lies in characterizing a set of assumptions that allow a 
reduced-order saddle-point representation of the 
optimal policy. These include a realizability assumption and a newly identified \emph{coherence 
assumption} about the subspaces used for approximation. Our main positive result is showing 
that these conditions are sufficient for constructing an algorithm that outputs an 
$\varepsilon$-optimal policy with runtime guarantees of $\wt{\OO}\pa{\taumix^2 N^3/ \varepsilon}$, 
where $N$ is the number of variables in the relaxed optimization problem, and $\taumix$ is a 
notion of mixing time. Our approach is based on the celebrated Mirror Prox algorithm of 
\citet{Nem04} (see also \citealp{Kor76}). We complement our positive results by showing that our 
newly defined coherence assumption is necessary for the relaxed saddle-point approach to be viable: 
we construct a simple example violating the assumption, where achieving full optimality on the 
relaxed problem leads to a suboptimal policy.

We are not the first to consider saddle-point methods for optimization in Markov decision 
processes. \citet{Wan17b} proposed variants of Mirror Descent to solve the original 
saddle-point problem without relaxations and provide runtime guarantees of 
$\wt{\OO}\bpa{\pa{\alpha \taumix}^2 |\X||\A| / \varepsilon^2}$, where $\X$ and $\A$ are the 
finite state and action spaces, and $\alpha$ is a parameter that 
characterizes the uniformity of the stationary distributions of every policy. Specifically, their 
assumption implies\footnote{The actual assumption made by \citet{Wan17b} is even more 
restrictive.} that for the stationary distribution $\stdist_\pi$ any policy $\pi$, one has 
$\frac{\max_x \stdist_\pi(x)}{\min_{x'} \stdist_\pi(x')} \le \alpha$. 
In most cases of practical interest, this ratio is at least as large as $|\X|$ (e.g., when there 
are states that some policies visit with constant probability), and can easily be exponentially 
large in $|\X|$, or even infinite if the underlying MDP has transient states. When specialized 
to this setting, our bounds replace $\alpha^2$ by the much more manageable $|\X|$ and also improve 
the dependence on $\varepsilon$ from $1/\varepsilon^2$ to $1/\varepsilon$. One downside of our 
method is that we need full access to the transition probabilities of the MDP, whereas the 
algorithm 
of \citet{Wan17b} only requires a generative model.

The linearly relaxed saddle-point problem we consider was first studied by \citet{LBS18} and
\citet{CLW18}. Our runtime guarantees improve over the ones claimed by \citet{CLW18} in a similar 
way as our first set of results improve over those of \citet{Wan17b}. Notably, their results 
still feature a factor of $\alpha^2$, which generally depends on the size of the original state 
space rather than the number of features, rendering these guarantees void of meaning in very large 
state spaces. In contrast, our bounds replace this factor by the number of features $N$. 
Furthermore, our characterization highlighting the importance of the coherence assumption discussed 
above hints at some potential technical issues with the results of \citet*{CLW18}, who claimed 
convergence to the optimal policy \emph{without the coherence assumption}.

The rest of the paper is organized as follows. After providing background 
on the saddle-point formulation of MDP optimization in Section~\ref{sec:prelim}, we describe the 
relaxed saddle-point problem in Section~\ref{sec:approx}.
Section~\ref{sec:alg} presents our algorithm and its performance guarantees, and 
Section~\ref{sec:analysis} provides a sketch of the 
proofs. We conclude by providing a simple numerical illustration of our method in 
Section~\ref{sec:experiments} and discuss our results in Section~\ref{sec:conc}.

\paragraph{Notation.} Inner products over vector spaces will be denoted by $\iprod{\cdot}{\cdot}$. 
We use $\Delta_{\S}$ to denote the set of probability distributions on the 
finite set $\S$: $\Delta_{\S} = \ev{p \in \real^\S_+: \sum_{s\in\S}p(s) = 1}$. Sums spanning over 
the spaces $x\in\X$ and $a\in\A$ will be simply denoted by $\sum_{x}$ or $\sum_a$.

\section{Preliminaries}\label{sec:prelim}
Consider an undiscounted Markov decision process  $M = (\mathcal{X}, 
\mathcal{A}, P, r)$, where  $\mathcal{X}$ is the finite state space, $\mathcal{A}$ is the finite 
action space, $P$ is the transition function with $P(x'|x,a)$ denoting the probability of 
moving to state $x'\in\mathcal{X}$  from state $x\in\mathcal{X}$ when taking action 
$a\in\mathcal{A}$ and $r$ is the reward function mapping state-action pairs to rewards 
with $r(x,a)$ denoting the reward of being in state $x$ and taking action $a$. We assume that 
$r(x,a)\in[0,1]$ for all $x,a$. In each round $t$, 
the learner observes state  $x_t\in\X$, selects action $a_t\in\mathcal{A}$, moves to the next state 
$x_{t+1} \sim P(\cdot|x_t,a_t)$, and obtains reward $r(x_t, a_t)$. 

In this paper we focus on the infinite-horizon average-reward scenario where the goal of the 
learner is to select its actions $a_t$ in a way that maximizes the average reward per time step, 
$\lim \inf_{t\rightarrow \infty} \mathbb{E}\left[\frac1{T}\sum_{t=1}^{T}r_t(x_t,a_t)\right]$. 
We will work with randomized stationary policies with $\pi(a|x)$ denoting the 
probability of taking action $a$ in state $x$. 
Under technical assumptions discussed 
shortly, each such policy $\pi$ generates a unique stationary state distribution
$\stdist_\pi\in\Delta_{\X}$ over the state space satisfying $\stdist_\pi(x) = \lim_{t\rightarrow 
\infty} 
\PP{x_t=x}$ for all $x$ when the trajectory $(x_t)_t$ is generated by following policy $\pi$. 
Similarly, each policy $\pi$ generates a stationary state-action distribution 
$\mu_\pi\in\Delta_{\X\times\A}$ satisfying $\mu_\pi(x,a) = \lim_{t\rightarrow \infty} 
\PP{x_t=x,a_t=a} = \stdist_\pi(x) \pi(a|x)$. Given these definitions, it can be easily shown that 
the 
average-reward of a policy $\pi$ can be written as
\[
 \rho_\pi = \lim \inf_{t\rightarrow 
\infty}\mathbb{E}_\pi\left[\frac1{T}\sum_{t=1}^{T}r_t(x_t,a_t)\right] = \sum_{x,a} \mu(x,a) 
r(x,a),
\]
where the notation $\EEpi{\cdot}$ indicates that the trajectory $(x_t,a_t)_t$ was generated by 
following policy $\pi$: $a_t\sim\pi(\cdot|x_t)$ and $x_{t+1} \sim P(\cdot|x_t,a_t)$. Under our 
assumptions, the optimal policy can be shown to be a stationary one; we will denote its average 
reward as $\rho^* = \max_\pi \rho_\pi$.
Thus, one can show that finding the optimal policy is equivalent to solving the following linear 
program:
\begin{eqnarray*}
 & \text{maximize} &\sum_{x,a} \mu(x,a) r(x,a)\\
 & \text{s.t.} &\mu \in \Delta_{\X\times\A}, \qquad \sum_{a'} \mu(x',a') = \sum_{x,a} 
P(x'|x,a)\mu(x,a) \quad(\forall x'\in\X).
\end{eqnarray*}
To simplify our notation, we will represent $\mu$ and $r$ by 
$\left|\X\times\A\right|$-dimensional vectors and also define the 
$\left|\X\times\A\right|\times|\X|$-dimensional matrix $Q$ with entries 
$Q_{(x,a),x'} = P(x'|x,a) - \II{x'=x}$.
Then, one can easily see\footnote{This can be seen, e.g., by introducing the KKT multipliers for 
the constraints in the linear program.} that
solving the linear program stated above is equivalent to finding the following \emph{saddle point}:
\begin{equation}\label{eq:saddle}
 \min_{v\in\real^{|\X|}} \max_{\mu\in\Delta} \mathcal{L}(v,\mu)=\min_{v\in\real^{|\X|}} 
\max_{\mu\in\Delta} \iprod{\mu}{Q v}+\iprod{\mu}{r}.
\end{equation}
Here, we introduced the \emph{Lagrangian function} $\LL$ and the shorthand $\Delta = 
\Delta_{\X\times\A}$. Optimal solutions $(v^*,\mu^*)$ to the above saddle-point problem are easily 
seen to correspond to 
the stationary distribution $\mu^*$ of the optimal policy and the \emph{optimal differential value 
function} $v^*$ (also known as the optimal bias function, cf.~\citealp{Puterman1994}). Besides the 
full saddle-point optimization problem, we will consider a relaxed version based on the 
introduction feature maps. Details on this variant are provided in Section~\ref{sec:approx}.

We will make two structural assumptions about the underlying Markov decision process. The first of 
these guarantees the existence of stationary distributions for all policies. 
\begin{assumption}[Uniform ergodicity]\label{ass:mixing}
Every policy $\pi$ generates an ergodic Markov chain. Specifically, letting $P_\pi$ be the 
transition operator of $\pi$ defined as the matrix with elements $P_\pi(x'|x) = \sum_a \pi(a|x) 
P(x'|x,a)$, and $\stdist,\stdist'$ be any two distributions over $\X$, the following 
inequality is satisfied for some $C,\tau > 0$ and for all $k$:
\[
 \onenorm{\pa{\stdist - \stdist'}P_\pi^k} \le C e^{-k/\tau}\onenorm{\stdist - \stdist'}.
\]
\end{assumption}
We say that our MDP is \emph{uniformly ergodic} if it satisfies Assumption~\ref{ass:mixing}.
Notice that this assumption is significantly weaker than the $1$-step mixing assumption 
often made in the related literature \citep{even-dar09OnlineMDP,neu14o-mdp-full}. It is easily 
shown to hold when all policies induce aperiodic and irreducible Markov chains---see 
Theorem~4.9 in \citet{LPW17} for a proof. Clearly, this assumption immediately implies that every 
policy admits a unique stationary distribution as required in the discussion above. In what follows 
below, we will often use the notation $\taumix = 2C\pa{\tau+1}$ and refer to this quantity as the 
\emph{mixing time} of the MDP\footnote{Note that this is just one of many possible definitions of 
a mixing time, see, e.g., \citet{Sen2006,LPW17}.}.

Given this assumption and the above definitions, we can establish a number of useful facts about 
the optimal solutions $(v^*,\mu^*)$ to the saddle-point problem~\eqref{eq:saddle}. We first note 
that an optimal policy $\pi^*$ can be extracted from $\mu^*$ in the states where $\mu^*(x,\cdot)>0$ 
as $\pi^*(a|x) = \frac{\mu^*(x,a)}{\sum_{a'} \mu^*(x,a')}$. 
Regarding $v^*$, the following proposition summarizes some of its most important 
properties:
\begin{prop}\label{prop:vprops}
Let $(v^*,\mu^*)$ be a solution of the problem~\eqref{eq:saddle}. Then, $v^*$ satisfies the 
following properties:
\begin{itemize}
 \item $v^*$ satisfies the Bellman optimality equations $v^*(x) = r(x) - \rho^* 
+ \sum_{x'}P(x'|x,a) v^*(x')$ for all $x$; for any $c\in\real$, $v^* + c$ is also a solution 
to~\eqref{eq:saddle};
 \item for any $x,x'$, $|v^*(x) - v^*(x')| \le \taumix = 2C\pa{\tau+1}$.
\end{itemize}
\end{prop}
All of these properties can be proven by standard arguments; we refer the reader to Lemma~1 in 
\citet{Wan17b} for a proof of the first item and Lemma~3 in \citet{neu14o-mdp-full} for a proof of 
the second one.

\section{The linearly relaxed saddle-point problem}\label{sec:approx}
While one can directly derive optimization algorithms to solve the saddle-point 
problem~\eqref{eq:saddle}, such a direct approach would suffer from serious scalability issues due 
to the sheer number of variables involved in the problem: the size of the objects of interest $\mu$ 
and $v$ are linear in the size of the state space, which results in prohibitive memory and 
computation costs for most algorithms. To address this issue, we study a \emph{linearly relaxed} 
version of the full saddle-point problem that reduces the order of the original 
optimization problem by linearly parametrizing the variables $v$ and $\mu$ through two sets of 
\emph{feature maps}. 
Formally, we consider the matrices $F$ of size $|\X| \times N$ and $W$ of size 
$M\times|\X\times\A|$, introduce the new optimization variables $y\in\real^M$ and $u\in 
\real^N$, and use these to (hopefully) approximate the solutions to~\eqref{eq:saddle} as 
$\mu^*\approx yW$ and $v^*\approx Fu$.
For a tractable problem formulation, we will assume that the rows of $W$ are non-negative and sum 
to one: $W_{m,x}\ge 0$ for all $x,m$ and $\sum_x W_{m,x} = 1$ for all $m$. We will also assume that 
all entries of $F$ are bounded by $1$ in absolute value. These conditions enable us to optimize $y$ 
over the probability simplex $\tD = \Delta_{[M]}$ and to formulate our relaxed saddle-point problem 
as 
\begin{equation}\label{eq:saddle_relax}
 \min_{u\in\real^{N}} \max_{y\in\tD} \wt{\mathcal{L}}(u,y)= \min_{u\in\real^{N}} \max_{y\in\tD} 
 \iprod{W\transpose y}{ Q Fu}+ \iprod{W\transpose y}{ r}.
\end{equation}
The relaxed optimization problem above has been studied before by \citet{LB15,LBS18}, and 
\citet{CLW18}. \citet{LB15,LBS18} studied the relaxed linear program 
underlying~\eqref{eq:saddle_relax} as a natural extension of the classic relaxed LP analyzed by 
\citet{FR03}, and have focused on understanding the discrepancies between the optimal value 
function and the relaxed value function attaining the minimum in the above expression. On the other 
hand, \citet{CLW18} focused on proposing stochastic optimization algorithms and analyzing the rate 
of convergence to the optimum, but provide little insight about the quality of the optimal solution 
of the relaxed problem. 

One of our main goals in the present paper is to obtain a better understanding of the effects of 
approximation on the policies that can be obtained through approximately solving the the relaxed 
saddle-point problem~\eqref{eq:saddle_relax}. 
One peculiar challenge associated with our setting is that it is not enough to ensure that the 
values of $\wt{\mathcal{L}}$ and $\LL$ are close at their 
respective saddle points, but we rather need to understand the performance of the policy
extracted from the optimal solution $y^*$. Precisely, defining the policy extracted from $y$ as
\[
 \pi_y(a|x) = \frac{(W\transpose y)(x,a)}{\sum_{a'}(W\transpose y)(x,a')}
\]
for all $x,a$, and the corresponding stationary distribution as $\mu_y$ induced in the original 
MDP, 
we are interested in the suboptimality gap $\iprod{\mu^* - \mu_{y^*}}{r}$.
In the present paper, we focus on identifying assumptions on the feature maps that allow the 
computation of true optimal policies with (almost) zero suboptimality gap. Specifically, we will 
show that the following two assumptions have a decisive role in making this gap small:
\begin{assumption}[Realizability]\label{ass:real}
 The optimal solution is realizable by the feature maps: there exists $\pa{u^*,y^*}$ such that $v^* 
= Fu^*$ and $\mu^* = W\transpose y^*$. Additionally, $\infnorm{u^*} \le U\taumix$ holds for 
some $U>0$.
\end{assumption}
\begin{assumption}[Coherence]\label{ass:coh}
 The image of the set $\tD$ under the map $Q\transpose W\transpose$ is included the column space of 
$F$: for all $y\in\tD$ such that $Q\transpose W\transpose y \neq 0$, there exists a $u\in\real^N$ 
such that $\iprod{Q\transpose W\transpose y }{F u} \neq 0$. Additionally, for all 
$v\in\real^{|\X|}$ with $\infnorm{v}\le 1$, there exists a $u\in\real^N$ with $\infnorm{u}\le 
U$ such that $\iprod{Q\transpose W\transpose y}{Fu} = \iprod{Q\transpose W\transpose y}{v}$.
\end{assumption}
The second condition of each assumption is to ensure that the columns of $F$ are 
well-conditioned and are satisfied if the columns form an orthonormal basis. While 
realizability may already seem sufficient for the relaxed problem to be a 
good enough approximation of the original one, we argue that the second assumption is also 
necessary for the relaxation scheme to be reliable. Specifically, the following theorem shows that 
in the absence of the coherence assumption, near-optimal solutions to the relaxed saddle-point 
problem~\eqref{eq:saddle_relax} can still lead to suboptimal policies in the original MDP. 
\begin{theorem}
 For any $\varepsilon>0$, there exists an MDP with relaxations $W,F$ satisfying 
Assumption~\ref{ass:real} and violating Assumption~\ref{ass:coh}, and a solution $(\wh{u},\hye)$ 
simultaneously satisfying
\[
\LL(F\wh{u},\mu^*) - \LL(v^*,W\transpose\hye) = \varepsilon
\]
and 
\[
 \iprod{\mu^* - \mu_{\hye}}{r} = 2/3.
\]
\end{theorem}

\begin{wrapfigure}[14]{r}{0.35\textwidth}
\includegraphics[width=.35\textwidth]{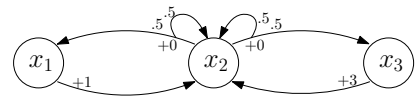} 
\caption{Three-state MDP for illustrating the necessity of the coherence assumption. Transitions 
from $x_2$ are stochastic with probability $1/2$ of staying in $x_2$ and moving to $x_1$ and $x_3$ 
otherwise, depending on the chosen action. All other transitions are deterministic. Rewards are 
given as a function of the state as $r(x_1) = 1$, $r(x_2) = 0$ and $r(x_3) = 3$.}\label{fig:cex}
\end{wrapfigure}
\begin{proof}
The proof is based on constructing an MDP with three states $x_1$ (left), $x_2$ (middle) and $x_3$ 
(right) and two actions $a_l$ and $a_r$ corresponding to moving ``left'' or ``right'', 
respectively. The transition probabilities and rewards are as shown on Figure~\ref{fig:cex}. It is 
easy to see that the optimal policy is to take action $a_r$ in state $x_2$, which yields the 
optimal stationary state-action distribution 
\begin{align*}
\mu^*=& 
\pa{\mu(x_1,a_r),\mu(x_2,a_l),\mu(x_2,a_r),\mu(x_3,a_l)}\transpose
=\pa{0, 0, \frac 13,\frac 23}\transpose
\end{align*}
and the optimal average reward $\rho^* = 1$. The optimal value function can be shown to be 
$v^* = \pa{-1,-1,1}\transpose$. For the relaxation, define $F = v^*$ and $W$ as the identity map 
so that the realizability assumption is clearly fulfilled with $y^* = \mu^*$ and $u^* = 1$. Now, 
choosing $\wh{y} = \pa{1, 0, 0, 0}\transpose$ results in
\[
 \iprod{W\transpose\wh{y}}{QFu} = \begin{pmatrix} 1&0&0&0 \end{pmatrix} \begin{pmatrix} 
-1 & 1 & 0 \\
 1/2 & -1/2 & 0 \\
0 & -1/2 & 1/2 \\
 0 & 1 & -1 
\end{pmatrix} \begin{pmatrix} 
-1 \\
-1 \\
1\end{pmatrix} u 
= \begin{pmatrix} 1&0&0&0 \end{pmatrix} 
\begin{pmatrix} 
0 \\
0 \\
1 \\
2 
\end{pmatrix} u
= 0\cdot u
\]
for any $u$. Observing that taking $v = (-1,1,0)\transpose$ gives 
$\iprod{W\transpose\wh{y}}{Qv} = 2$, we see that the coherence assumption is violated since there 
exists no $u$ such that the condition $\iprod{W\transpose\wh{y}}{Qv} 
= \iprod{W\transpose\wh{y}}{QFu}$ is satisfied. Furthermore, it is easy to see that $(\wh{y},u)$ 
for any $u$ is an optimal solution to~\eqref{eq:saddle_relax} with value $\rho^* = 1$ since
\begin{align*}
    \wt{\mathcal{L}}(u,\wh{y})&= \wh{y}\transpose WQFu + \wh{y}\transpose W r   = \begin{pmatrix} 
1&0&0&0 \end{pmatrix}\begin{pmatrix} 
1 \\
0 \\
0 \\
3 \\
\end{pmatrix}
= 1.
\end{align*}
showing that $(\wh{y},u)$ with any $u$ is also an optimal solution to the relaxed saddle-point 
problem~\eqref{eq:saddle_relax}. The resulting optimal state-action distribution $\wh{\mu} = 
\wh{y}W = \wh{y}$ is clearly not a stationary distribution.

To conclude the proof, fix any $\varepsilon$ and consider $\hye = \pa{1 - 
\varepsilon, \varepsilon, 0, 0}\transpose$ and any $\wh{u}$. Noticing that 
$\iprod{W\transpose\hye}{QFu} = 0$ holds for all $u$, the duality gap 
associated with $(\wh{u},\hye)$ can be seen to be 
\begin{align*}
 \LL(F\wh{u},\mu^*) - \LL(v^*,W\transpose\hye) &= \begin{pmatrix} 
0&0&2/3&1/3 \end{pmatrix}\begin{pmatrix} 
1 \\
0 \\
0 \\
3 \\
\end{pmatrix} 
-  
\begin{pmatrix} 
1-\varepsilon&\varepsilon&0&0 
\end{pmatrix}
\begin{pmatrix} 
1 \\
0 \\
0 \\
3 \\
\end{pmatrix}
= 1 - (1-\varepsilon) = \varepsilon.
\end{align*}
The policy $\pi_{\hye}$ extracted from the state-action distribution $\hye$ takes action $a_l$ in 
state $x_2$, which results in an average reward of $2/3$. These two statements together prove the 
theorem.
\end{proof}

\section{Algorithm and main results}\label{sec:alg}
In this section, we provide our main positive results: deriving strong performance guarantees for 
policies derived from approximate solutions  of~\eqref{eq:saddle_relax} under 
Assumptions~\ref{ass:real} and~\ref{ass:coh}. Our algorithm attaining these guarantees is 
based on the Optimistic Mirror Descent framework proposed by \citet{RS13,RS13b}, and more 
specifically on its variant known as Mirror Prox due to \citet{Nem04} (see also 
Sections~4.5 and~5.2.3 in \citet{Bub15} for an easily accessible overview of this method). 

For a generic description of Mirror Prox on a convex set $\Z$, we let $G:\Z\ra\real$ be a 
monotone operator 
satisfying $\iprod{G(z)-G(z')}{z-z'}\geq0$ for all $z,z'\in\Z$, and let $\Phi:\Z\ra\real$ be 
a $\sigma$-strongly convex regularization function under some norm $\norm{\cdot}$ with its 
corresponding Bergman divergence $\DD{z}{z'}= \Phi(z)-\Phi(z')-\iprod{\nabla \Phi(z')}{z-z'}$. 
Mirror Prox computes a sequence of iterates with $z_1\in \arg \min \Phi(z)$ and
\begin{equation} \label{eq:MP}%\label{opt}
\begin{split}
     &\wh{z}_{t+1}=\arg \min_\mathcal{Z} \eta \iprod{G(z_t)}{z} + D_\Phi(z,z_t)\\
    &z_{t+1}=\arg \min_\mathcal{Z} \eta \iprod{G(\wh{z}_{t+1})}{z} + D_\Phi(z,z_t).
\end{split}
\end{equation}
The first of these steps is often referred to as an \emph{extrapolation step}. A simpler version of 
this algorithm not involving such an extrapolation step is commonly known as Mirror Descent 
\citep{NY83,BT03,Bub15}.
This step serves to enhance the stability of the algorithm, and indeed Mirror Prox can be shown to 
enjoy favorable convergence properties in the problem setting described above.

We instantiate the Mirror Prox method to address the relaxed saddle-point problem as follows. Our 
optimization variables will be $z=(u,y)$ and the monotone operator $G$ will be chosen as 
\begin{equation}\label{eq:Gdef}
G(z)= \left( \begin{array}{c}
\nabla_v \wt{\LL}\\
-\nabla_\mu \wt{\LL} \end{array} \right)
=\left( \begin{array}{c}
F\transpose Q\transpose W\transpose y \\
-Wr-WQFu \end{array} \right).
\end{equation}
We will use the regularization function 
\[
 \Phi(z) = \frac{1}{2} \norm{u}_2^2 + \sum_{j=1}^M y_j \log y_j,
\]
that is, a linear combination of the squared $2$-norm of the value-function parameters $u$ and the 
Shannon entropy of the distribution $y$. Clearly, $\Phi$ is $1$-strongly convex on $\Z$ with 
respect 
to the norm $\norm{z}^2 = \norm{u}_2^2 + \norm{y}_1^2$.
Given the above specifications, the updates of our algorithm can be written as
\begin{eqnarray}\label{eq:update1}
    &\wh{u}_{t+1} =   u_t - \eta F\transpose Q\transpose W\transpose y_t,
    &\wh{y}_{t+1,i} \propto y_{t,i} e^{\eta \pa{(Wr)_i + \pa{WQFu_t}_i}}\\
    &u_{t+1} =   u_t - \eta F\transpose Q\transpose W\transpose \wt{y}_{t+1},
    &y_{t+1,i} \propto y_{t,i} e^{\eta \pa{(Wr)_i + \pa{WQF\wh{u}_{t+1}}_i}},\label{eq:update2}
\end{eqnarray}
where we used the notation ``$\propto$'' to signify that $\wh{y}_{t+1}$ and $y_{t+1}$ are 
normalized multiplicatively after each update so that $\sum_j y_{t+1,j} = 1$ is satisfied.
Also introducing the notations $\overline{y}_T = \frac{1}{T} \sum_{t=1}^T y_t$ and $\overline{u}_T 
= \frac{1}{T} \sum_{t=1}^T \wh{u}_t$, the algorithm outputs the policy extracted from the 
distribution $\overline{y}_T$: $\pi_T = \pi_{\overline{y}_T}$. 
Letting $\stdist_T = \stdist_{\pi_T}$ be the stationary distribution associated with $\pi_T$, the 
corresponding average reward can be written as $\rho_T = \sum_{x,a} \stdist_T(x) \pi_T(a|x) 
r(x,a)$. 
The following theorem presents our main result regarding the suboptimality of the resulting policy 
in terms of its average reward.
\begin{theorem}\label{thm:mainapp}
Suppose that Assumptions~\ref{ass:mixing},~\ref{ass:real} and~\ref{ass:coh} hold and $\eta \le 
 1/4N$. Then, the average reward $\rho_T$ output by the algorithm satisfies
 \[
  \rho^* - \rho_T \le  \frac{11\taumix^2 U^2 N + 7\log M}{\eta T}.
 \]
 In particular, setting $\eta = 1/4N$, the bound becomes $\rho^* - \rho_T = 
\OO\pa{\frac{\taumix^2N^2U^2}{T}}$.
\end{theorem}
We note that this result can be tightened by a factor of $N$ if we further assume that the rows of 
$F$ are chosen as probability distributions. In the special case where $F$ and $W$ are the identity 
maps, the relaxed saddle-point problem becomes the original problem~\eqref{eq:saddle}, and our 
Assumptions~\ref{ass:real} and~\ref{ass:coh} are clearly satisfied with $U=1$. In this case, our
algorithm satisfies the following bound:
\begin{corollary}\label{thm:main}
 Suppose that Assumption~\ref{ass:mixing} holds, $W$ and $F$ are the identity maps, and $\eta \le 
1/4$. Then, the average reward $\rho_T$ of the policy output by our algorithm satisfies
 \[
  \rho^* - \rho_T \le  \frac{11\taumix^2 |\X| + 7\log\pa{|\X||\A|}}{\eta T}.
 \]
 In particular, setting $\eta = 1/4$, the bound becomes $
\rho^* - \rho_T = \wt{\OO}\pa{\frac{\taumix^2 |\X|}{T}}$.
\end{corollary}

A brief inspection of Equations~\eqref{eq:update1}-\eqref{eq:update2} suggests that each update of 
our algorithm can be computed in $\OO\pa{MN}$ time, the most expensive operation being computing 
the matrix-vector products $WQFu$ and $y\transpose WQF$. While this suggests that the algorithm may 
have runtime and memory complexity independent of the size of the state space, we note that 
exact computation of the matrix $WQF$ can still take $\OO\pa{|\X|^2|\A|}$ time in the worst 
case. This can be improved to $\OO\pa{K}$ when assuming that only $K$ entries of the transition 
matrix $P$ are nonzero, which can be of order $|\X||\A|$ in many interesting problems where the 
support of $P(\cdot|x,a)$ is of size $\OO\pa{1}$ for all $x,a$. We stress however that the matrix 
$WQF$ only needs to be computed \emph{once} as an initialization step of our algorithm. In 
contrast, a general algorithm like value iteration needs at least $\Theta\pa{K} = 
\Theta\pa{|\X||\A|}$ 
for computing \emph{each update}, showing a clear computational advantage of our method. Further 
discussion of computational issues is deferred to Section~\ref{sec:conc}.

\section{Analysis}\label{sec:analysis}
This section provides an outline of the analysis of our algorithm. At a high level, our analysis 
builds on some well-known results regarding the performance of Mirror Prox, including a classical 
bound on the \emph{duality gap} of the obtained solutions. The crucial challenge posed by our 
setting is connecting the duality gap on the saddle-point problem to a suboptimality gap of the 
extracted policies. The key innovation in our analysis is providing a new technique to connect 
these quentities through exploiting further properties of Mirror Prox. In what follows, we first 
provide some general tools that will be helpful throughout the proofs, and then provide the proof 
outline for Theorem~\ref{thm:mainapp}. 
Full proofs are provided in Appendix~\ref{app:proofs}. 

A central piece of our our analysis is the following useful lemma regarding the iterates computed 
by Mirror Prox:
\begin{lemma}\label{lem:OMD}
Let $\Phi$ be $\sigma$-strongly convex and $F$ be $L$-Lipschitz. Then, for all $t$, Mirror Prox 
guarantees
\[
 \eta \iprod{\hz_{t+1} - z{}}{G(\wh{z}_{t+1})}\le \DD{z{}}{z_{t}} - \DD{z{}}{z_{t+1}} - 
\frac{\sigma - 
\eta L}{4}\norm{z_{t+1} - z_t}^2.
\]
holds for every $z\in\Z$ and $t>0$.
\end{lemma}
The proof is based on standard arguments, see, for instance, Lemma~1 of \citet{RS13b}. We include 
it in Appendix~\ref{app:OMD} for completeness. This lemma has two important corollaries that we 
will crucially use throughout the analysis. The first one shows that the iterates remain bounded 
during the optimization procedure.
\begin{corollary}\label{cor:normbound}
 Let $z^* = \pa{u^*,y^*}$ be any solution to $\max_{y}\min_{u} \wt{\LL}\pa{u,y}$ and suppose that 
the conditions of Lemma~\ref{lem:OMD} hold. Then, for all $t$, Mirror Prox guarantees 
\begin{align*}
\DD{z^*}{z_{t}} \le \DD{z^*}{z_0}. 
\end{align*}
\end{corollary}
The proof follows from noticing that $z^*$, being an optimal solution to the saddle-point 
problem, satisfies the variational inequality $\iprod{\hz_{t+1} - z^*}{G(\wh{z}_{t+1})} \ge 0$.
The second corollary establishes a bound on the \emph{duality gap} evaluated at $(\bu_T,\by_T)$:
\begin{corollary}\label{cor:gapbound}
 Let $z = \pa{u,y}\in\Z$ be arbitrary and assume that $\eta \le \frac{\sigma}{2L}$. Then, 
Mirror Prox guarantees the following bound on the duality gap:
\begin{align*}
\mathcal{L}\left( \bu_T,y \right) - \mathcal{L}\left(u, \by_T \right) 
     \leq& \frac{D_\Phi(z, z_0)}{\eta T}.
\end{align*}
\end{corollary}
The proof easily follows by noticing that $\iprod{\hz_{t+1} - z{}}{G(\wh{z}_{t+1})}$ equals the 
duality gap evaluated at $(\hu_{t+1},\hy_{t+1})$, and summing the bound given in 
Lemma~\ref{lem:OMD}.

In order to apply the above tools to our problem, we first need to confirm that our objective is 
indeed smooth with respect to the norm $\norm{z}^2 = \norm{u}_2^2 + \norm{y}_1^2$. The following 
lemma establishes this property.
\begin{lemma}\label{lem:smoothness_relax}
Let $K = \max_x \onenorm{F_{x,\cdot}}$. Then, the function $\wt{\LL}$ is $2K$-smooth with 
respect to $\norm{\cdot}$. 
\end{lemma}
The proof is provided in Appendix~\ref{app:smoothness_relax}. Notably, this lemma implies that the 
$\wt{\LL}$ is $2$-smooth when the rows of $F$ form probability distributions. In the worst case, 
however, when we only assume that the entries of $F$ are bounded in absolute value by $1$, the 
smoothness constant can be as large as $2N$. In what follows, we will assume that $\eta \le 
1/(4K)$.

We proceed by appealing to the realizability assumption to choose $x=(u^*,y^*)$ such that 
$z=(v^*,\mu^*)=(Fu^*,W\transpose y^*)$, and observe that
\[
 \wt{\LL}\left( \bu_T,y^* \right) - \wt{\LL}\left(u^*, \by_T\right) 
 = \iprod{\mu^*}{QF\bu_T + r} - \iprod{W\transpose \by_T}{Qv^* + r} \le \frac{D_\Phi(z^*, 
z_0)}{\eta 
T}
\]
holds by virtue of Corollary~\ref{cor:gapbound} and the choice of $\eta$. Observing that 
$Q\transpose \mu^* = 0$ holds due to the stationarity of $\mu^*$ and reordering gives
\begin{equation}\label{eq:almostbound}
 \iprod{\mu^* - W\transpose\by_T}{r} \le \frac{D_\Phi(z^*, z_0)}{\eta T} + 
\iprod{Q\transpose W\transpose\bmu_T}{v^*}.
\end{equation}
The remaining key question is how to relate $\iprod{W\transpose\by_T}{r}$ to the true average 
reward $\rho_T$ associated with the extracted policy. This is done with the help of the following 
lemma, one of our key results:
\begin{lemma}\label{lem:nubound}
Suppose that Assumption~\ref{ass:mixing} holds. Let $\mu$ be an arbitrary distribution over 
$\X\times\A$ and let $\pi_\mu$ be the policy extracted from $\mu$. Then, the average reward 
$\rho_\mu$ of $\pi_\mu$ satisfies
$\iprod{\mu}{r} - \rho_\mu \le \taumix\onenorm{Q\transpose\mu}$.
\end{lemma}
The proof is provided in Appendix~\ref{app:nubound}. Combining this result with the bound of 
Equation~\ref{eq:almostbound} and using that $\infnorm{v^*}\le \taumix$, we obtain
\begin{equation}\label{eq:rhobound}
 \rho^* - \rho_T \le \frac{D_\Phi(z^*, z_0)}{\eta T} + 2\taumix \onenorm{Q\transpose 
W\transpose\by_T}.
\end{equation}
Thus, it only remains to bound $\onenorm{Q\transpose W\transpose\by_T}$. In order to do this, we 
crucially use Assumption~\ref{ass:coh} that guarantees the coherence of the feature maps to prove 
the following result:
\begin{lemma}\label{lem:yWQbound}
Suppose that Assumptions~\ref{ass:real} and~\ref{ass:coh} hold. Then,
 \[
  \taumix \onenorm{ Q\transpose W\transpose \bar{y}_T} \le \frac{5\taumix^2 U^2 N + 
3 \log M}{\eta T}
 \]
\end{lemma}
The proof of this lemma is provided in Appendix~\ref{app:yWQbound}.
Combining the bound of this lemma with Equation~\eqref{eq:rhobound} and using $\DD{z^*}{z_0} \le 
\taumix^2 N + \log\pa{M}$ concludes our proof of Theorem~\ref{thm:main}.

\section{Numerical illustration}\label{sec:experiments}
In this section, we provide empirical results on two simple environment in order to illustrate our 
theoretical results, and specifically compare the performance of our algorithm with that of Mirror 
Descent and the classic value iteration algorithm. 

In the first example, we consider a rectangular $s\times s$ gridworld with one nonzero reward placed in 
state $x_r$, so that $r(x,a)=\mathcal{I}_{x=x_r}$. Once the agent arrives to $x_r$, it is 
randomly teleported to any of the other states with equal probability. In any other state, the agent 
can decide to move to a neighboring cell in any direction. The attempt to move in the desired 
direction is successful with probability $p$, otherwise the agent moves in the opposite direction 
with probability $1-p$. If the agent is in an edge of the grid and it makes an step in the direction
of the edge, it appears in the opposite edge.

Figure \ref{regrets} shows some results on a grid of side $s=10$, in the case when no features are used, 
so we optimize over the whole state-action space. We observe that the 
convergence of Mirror Prox is much faster than that of Mirror Descent, and that the last iterate
of MP converges very quickly to the optimum, achieving it after \emph{finitely many iterations}. We 
also note that for higher values of $\eta$ than the ones found to be safe in our bounds 
(at most 1/4), the algorithm is still stable and can lead to faster convergence to the optimum.

\begin{figure}[h]
    \centering
    \begin{minipage}{0.5\textwidth}
        \centering
        \includegraphics[width=0.9\textwidth]{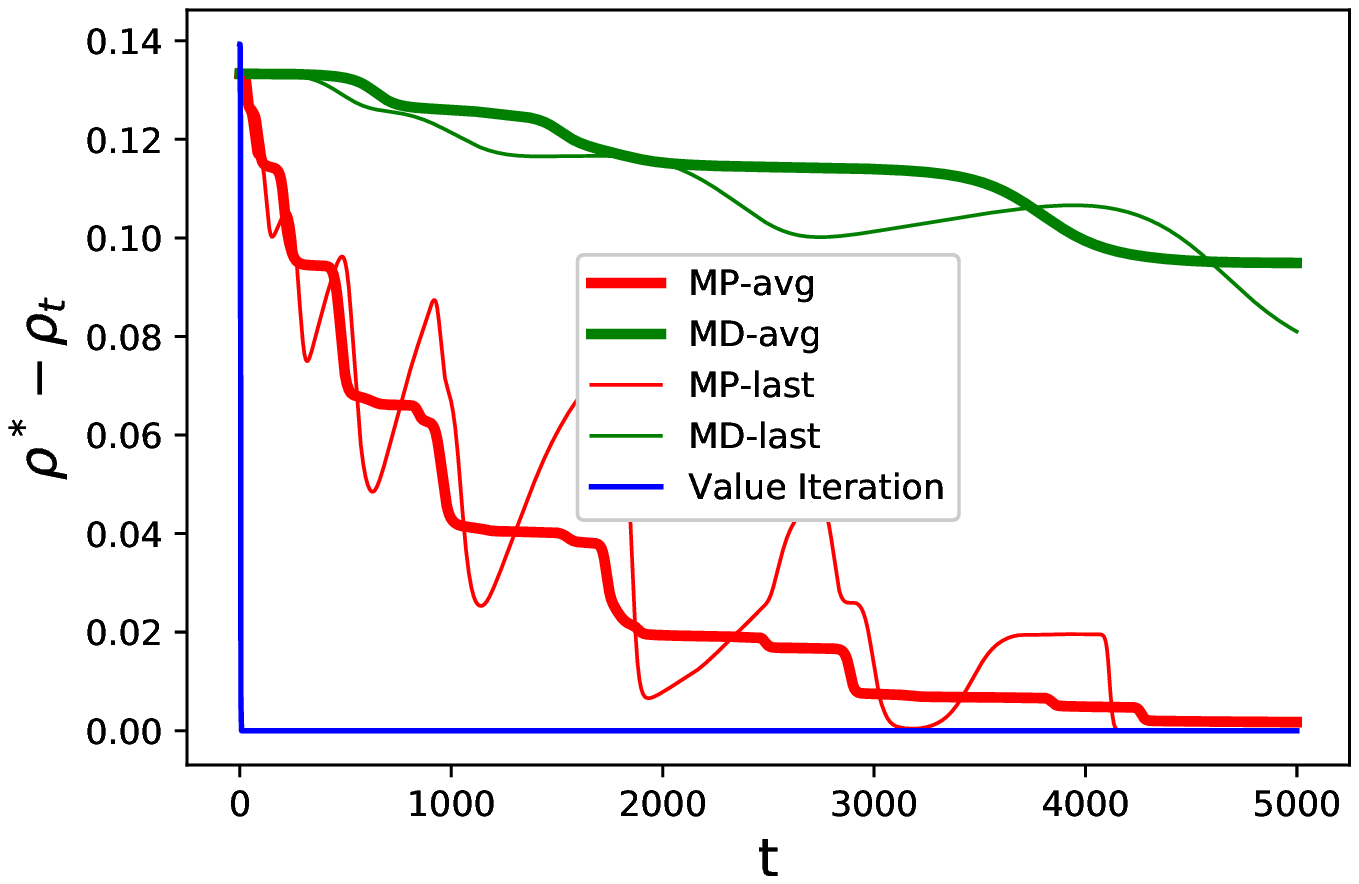} % first figure itself
        \subcaption{$p=0.9$, and $\eta=\frac 14$.}
    \end{minipage}\hfill
    \begin{minipage}{0.5\textwidth}
        \centering
        \includegraphics[width=0.9\textwidth]{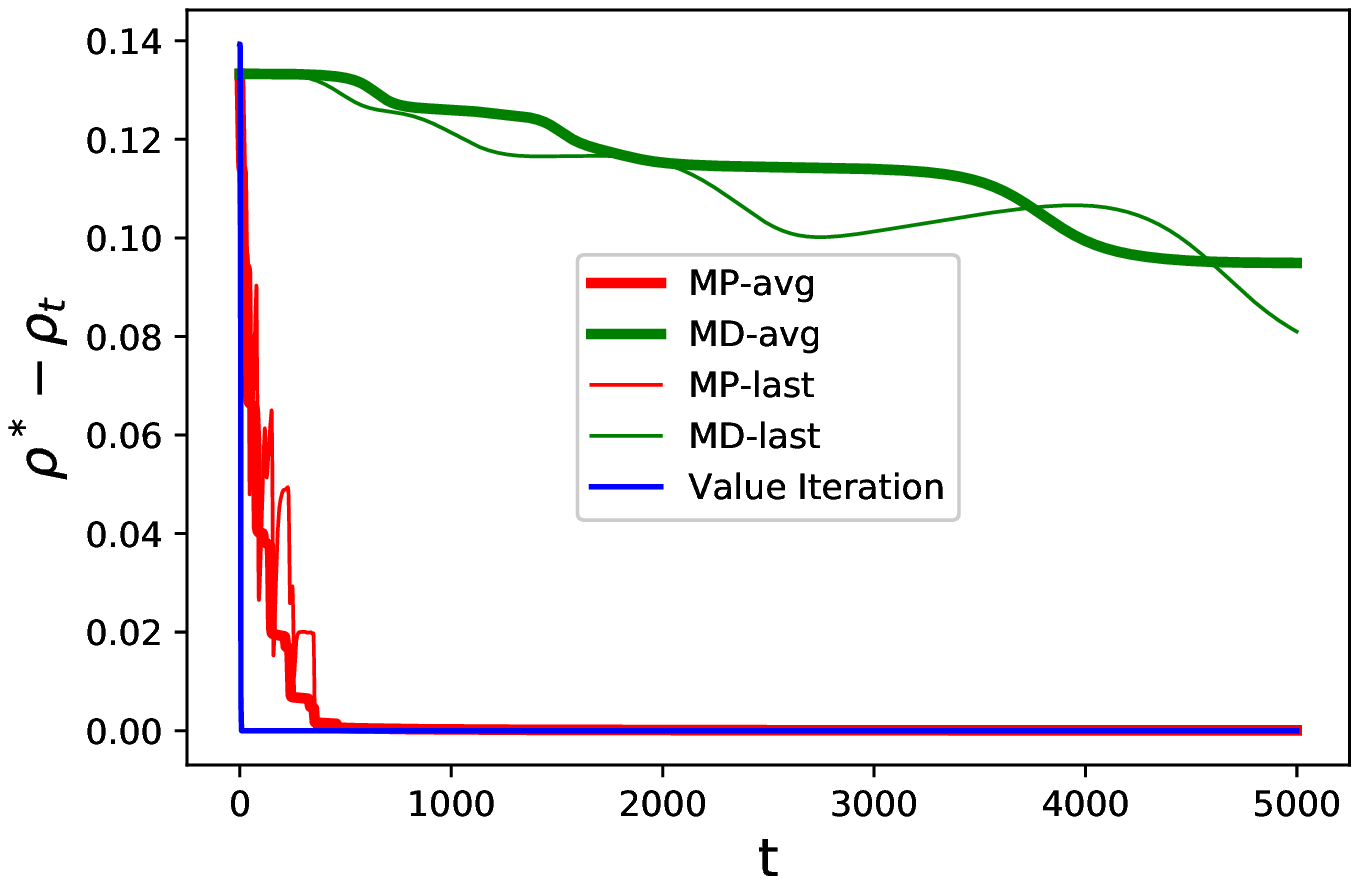} % second figure itself
        \subcaption{$p=0.9$, and $\eta=3$.}
    \end{minipage}
    \caption{Regret as a function of the number of iterations of MP, MD, and value 
iteration in a grid world example.}
\label{regrets}
\end{figure}

\begin{figure}[h]
\centering
\includegraphics[width=0.4\textwidth]{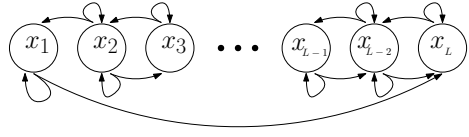} % first figure itself
\caption{Example of MDP.}
\label{featexample}
\end{figure}

In our second example, we show how the usage of good features can make MP converge faster than 
value iteration. 
We consider a sequence of states of length $L$ (see Figure~\ref{featexample}) with one nonzero 
reward placed in the first state so that $r_{(x, a)} = I_{x=x_1}L$. 
In states $x_2$ to $x_{N-1}$ the available actions are to go left and right, in state $x_1$ the only available action is to go to the 
last state ($x_L$), and in state $x_L$ the only available action is to go left. Each action has a 
probability $p$ of success and $1-p$ of remaining in the same state.

To test our algorithm in this environment, we built $W$ and $F$ taking advantage of the structure of 
the problem as follows: For $W$, we randomly generate a vector $c$ of length $L$ with entries 
being 1, 2 or 3. 
For $i\leq3$, $W^\top_{(x=i,a=\text{left}),j}$=1 if $c_j=i$ and 0 otherwise. After that we normalize the three rows, getting three homogeneous non-overlapping distributions. 
Doing this, we ensure that the realizability assumption is fulfilled for the $\mu_s$. We do the same 
for the ``\textit{right}'' action, and we add two more rows with 
random probability distributions over the whole set of state-action pairs. This makes for a total 
of 8 rows in $W$.

To build $F$,  we also randomly generate a vector $c$ of length $L$ with entries being 1, 2 or 3. 
For $i\leq 3$, $F_{j,i}=j$ if $c_j=j/L$ and 0 otherwise, to guarantee that the relaizability assumption is fulfilled for the $v_s$. We also add three random columns with 
random numbers between 0 and 1, in order to fulfill coherence with high probability. This results 
in a total of $5$ columns for $F$.

In Figure~\ref{regrets} we show the results obtained with value iteration and the linearly relaxed 
mirror prox, with $p=0.7$ and different lengths (10 and 100). 
While for value iteration the  number of iterations needed to converge is of the order of the number 
of states, it is independent of the size of the state space for our algorithm, and rather scales 
with the number of columns of the matrices $W$ and $F$. This simple example shows that with proper 
features, our algorithm  can actually beat value iteration, which by itself is not able to deal with 
features.

\begin{figure}[h]
    \centering
    \begin{minipage}{0.5\textwidth}
        \centering
        \includegraphics[width=0.9\textwidth]{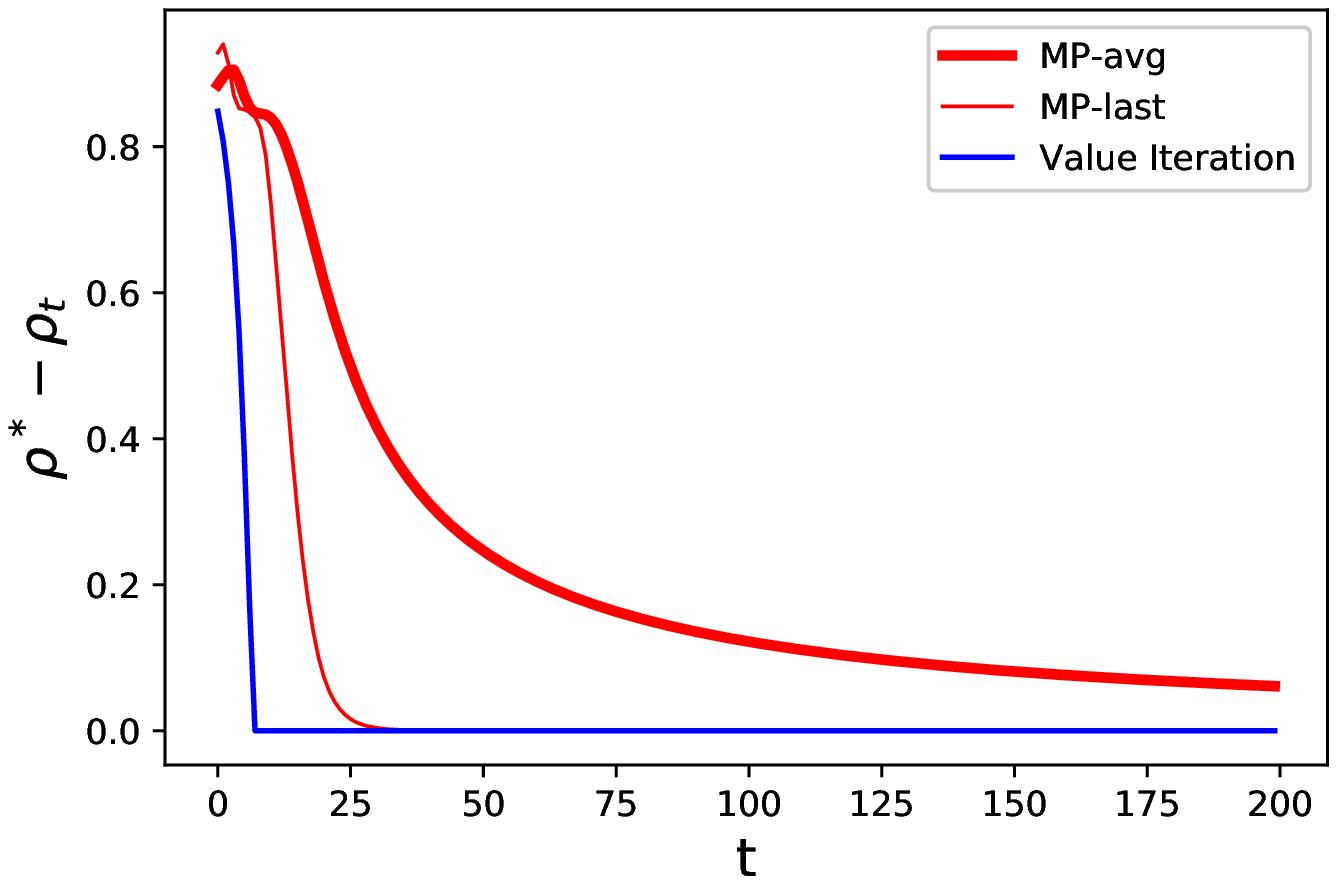} % first figure itself
        \subcaption{$p=0.7$, $\eta=\frac 14$ and $L=10$.}
    \end{minipage}\hfill
    \begin{minipage}{0.5\textwidth}
        \centering
        \includegraphics[width=0.9\textwidth]{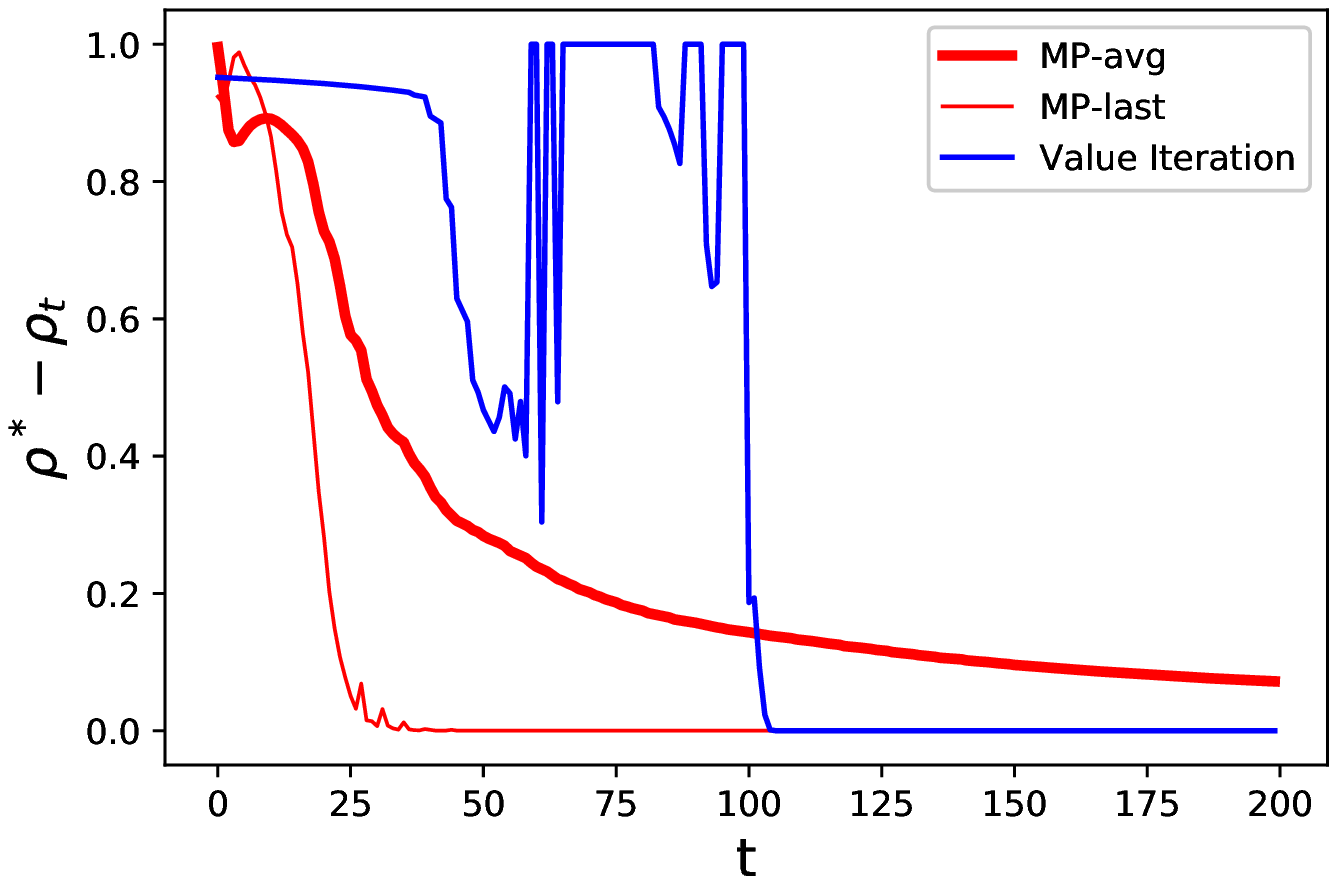} % second figure itself
        \subcaption{$p=0.7$, $\eta=\frac14$ and $L=100$.}
    \end{minipage}
    \caption{Suboptimality gap as a function of the number of iterations of MP and value 
iteration for $p=0.7$ and $\eta=0.25$}
\label{regrets}
\end{figure}

\section{Discussion}\label{sec:conc}
Our most important contributions concern the relaxed saddle-point problem~\eqref{eq:saddle_relax}, 
most notably including our discussion on the necessity and sufficience of the coherence assumption 
(Assumption~\ref{ass:coh}).
As we've mentioned earlier, several relaxation schemes similar to ours have been studied in the 
literature. In fact, relaxing the linear program underlying~\eqref{eq:saddle} through the 
introduction of the feature map $F$ for approximating the value function $v^*$ is one of the oldest 
ideas in approximate dynamic programming, originally introduced by \citet{SS85}. The effects of 
this approximation were studied by \citet{FR03} in the context of discounted Markov decision 
processes. A relaxation scheme involving both the feature maps $F$ and $W$ was considered by 
\citet{LB15,LBS18}. Both sets of authors carefully observed that introducing relaxations may make 
the linear program unbounded, and proposed algorithmic steps and structural assumptions of $F$ and 
$W$ to fight this issue. The results of these works are incomparable to ours since they focus on 
controlling the errors in approximating the optimal value function $v^*$ rather than controlling 
the 
suboptimality of the policies output by the algorithm. Interestingly, the widely popular REPS 
algorithm of \citet{peters10reps} is also originally derived from the relaxed linear program 
analyzed by \citet{FR03}, even if this connection has not been pointed out by the authors.

The work of \citet{CLW18} is very close to ours in spirit. \citeauthor{CLW18} consider a variation 
of the relaxed saddle-point problem~\eqref{eq:saddle_relax} with $W$ being block-diagonal 
with $F\transpose$ in each of its blocks, and claim convergence results for their algorithm to the 
optimal policy under only a realizability assumption. Unfortunately, their choice of $W$ does not 
necessarily ensure that the coherence assumption holds, which raises concerns regarding the 
generality of their guarantees. Indeed, the results of \citeauthor{CLW18} require an additional 
assumption that implies that $\frac{\max_x \stdist_\pi(x)}{\min_{x'} \stdist_\pi(x')}$ remains 
bounded by a constant for any policy $\pi$, which is extremely difficult to ensure in problems 
of practical interest. In fact, this ratio is already exponentially large in $|\X|$ 
in very simple problems like the one we consider in our experiments. 
Additionally, the analysis of \citeauthor{CLW18}  is based on the potentially erroneous claim 
that under the realizability assumption, the representation $(u^*,y^*)$ of the original optimal 
solution $(v^*,\mu^*) = (Fu^*,W\transpose y^*)$ always remains an optimal solution to the relaxed 
saddle-point problem. It is currently unclear if this claim is indeed true, or to what extent their 
condition regarding the boundedness of stationary distribution can be relaxed. 

In any case, we believe that our coherence assumption is more fundamental than the previously 
considered conditions, and it enables a much more transparent analysis of optimization algorithms 
addressing the relaxed saddle-point problem~\eqref{eq:saddle_relax}. Beyond this 
particular positive result, our work also cleans the slate for further theoretical work on 
approximate optimization in Markov decision processes. Indeed, the form of our coherence assumption 
naturally invites the question: can we compute good approximate solutions to the original problem 
when our assumptions are only satisified approximately? Similar questions are not without precedent 
in the reinforcement-learning literature. Translated to our notation, classical results concerning 
the performance of (least-squares) temporal difference learning algorithms imply that the 
approximation errors are controlled by the projection error of $QFu^* + r$ to the column space 
of $F$ \citep{TvR97,BB96,LGM10}. When using more general function classes to approximate $v^*$, 
\citet{MS08} show that the approximation errors are controlled by the \emph{inherent Bellman error} 
of the function class, which captures an approximation property related to our coherence condition. 
Whether or not we can generalize our techniques to construct provably efficient algorithms under 
such milder assumptions remains an exciting open problem that we leave open for future research.

% \bibliographystyle{plainnat}
% \bibliography{ngbib,confs}
% 
% 

\appendix
\section{Ommitted proofs}\label{app:proofs}
\subsection{The proof of Lemma~\ref{lem:OMD}}\label{app:OMD}
The proof will rely on repeatedly using the so-called three-points identity that can easily be 
shown to hold for all points $x,y,z\in\Z$:
\[
 \DD{x}{y} = \DD{x}{z} + \DD{z}{y} + \iprod{\nabla\Phi(y) - 
\nabla\Phi(z)}{z - x}.
\]
We first use it to show
\begin{align*}
 \DD{z{}}{z_{t+1}} &= \DD{z{}}{z_{t}} - \DD{z_{t+1}}{z_{t}} + \eta \iprod{z{} - 
z_{t+1}}{\nabla\Phi(z_{t+1}) - \nabla\Phi(z_t)}
\\
 &\le \DD{z{}}{z_{t}} - \DD{z_{t+1}}{z_{t}} + \eta \iprod{z{} - 
z_{t+1}}{G(\hz_{t+1})},
\end{align*}
where we also used the first-order optimality condition for $z_{t+1}$ in the second step:
\begin{align*}
 \iprod{\nabla\Phi(z_t) - \nabla\Phi(z_{t+1}) - \eta G(\hz_{t+1})}{z_{t+1} - z{}} \ge 0.
\end{align*}
Furthermore, we have
\begin{align*}
 \iprod{z{} - z_{t+1}}{G(\wh{z}_{t+1})} &= \iprod{z{} - \hz_{t+1}}{G(\wh{z}_{t+1})} + 
 \iprod{\hz_t - z_{t+1}}{G(\wh{z}_{t+1})}. 
\end{align*}
Using this bound together with the three-points identity
\[
 \DD{z_{t+1}}{z_t} = \DD{z_{t+1}}{\hz_{t+1}} + \DD{\hz_{t+1}}{z_t} + \iprod{\nabla\Phi(z_t) - 
\nabla\Phi(\hz_{t+1})}{\hz_{t+1} - z_{t+1}},
\]
we obtain
\begin{align*}
 \DD{z{}}{z_{t+1}} &\le \DD{z{}}{z_{t}} - \DD{z_{t+1}}{z_{t}} + \eta \iprod{\hz_{t+1} - 
z_{t+1}}{G(\wh{z}_{t+1})} + \eta \iprod{z{} - \hz_{t+1}}{G(\wh{z}_{t+1})}
\\
& = \DD{z{}}{z_{t}} - \DD{z_{t+1}}{\hz_{t+1}} - \DD{\hz_{t+1}}{z_t} + \eta \iprod{z{} - 
\hz_{t+1}}{G(\wh{z}_{t+1})}
\\
&\qquad\qquad + \iprod{\nabla\Phi(z_t) - 
\nabla\Phi(\hz_{t+1}) - \eta G(\hz_{t+1})}{z_{t+1} - \hz_{t+1}}
\\
& = \DD{z{}}{z_{t}} - \DD{z_{t+1}}{\hz_{t+1}} - \DD{\hz_{t+1}}{z_t} 
\\
&\qquad\qquad + \iprod{\nabla\Phi(z_t) - \nabla\Phi(\hz_{t+1}) - \eta G(z_t)}{z_{t+1} - \hz_{t+1}}
+ \eta \iprod{G(z_t) - G(\hz_{t+1})}{z_{t+1} - \hz_{t+1}}
\\
&\qquad\qquad + \eta \iprod{z{} - \hz_{t+1}}{G(\wh{z}_{t+1})}
\\
& \le \DD{z{}}{z_{t}} - \DD{z_{t+1}}{\hz_{t+1}} - \DD{\hz_{t+1}}{z_t} 
 + \eta \iprod{G(z_t) - G(\hz_{t+1})}{z_{t+1} - \hz_{t+1}}
\\
&\qquad\qquad + \eta \iprod{z{} - \hz_{t+1}}{G(\wh{z}_{t+1})} ,
\end{align*}
where the last step follows from the fact that  $\hz_{t+1}$ satisfies the first-order optimality 
condition
\begin{align*}
 \iprod{\nabla\Phi(z_t) - \nabla\Phi(\hz_{t+1}) - \eta G(z_t)}{z_{t+1} - \hz_t} \le 0.
\end{align*}
Now, using the $\sigma$-strong convexity of $\Phi$ and the $L$-Lipschitz continuity of $F$, we 
obtain
\begin{align*}
 \DD{z{}}{z_{t+1}} & \le \DD{z{}}{z_{t}} - \DD{z_{t+1}}{\hz_{t+1}} - \DD{\hz_{t+1}}{z_t} 
 + \eta \iprod{G(z_t) - G(\hz_{t+1})}{z_{t+1} - \hz_{t+1}}
 \\
 &\qquad\qquad + \eta \iprod{z{} - \hz_{t+1}}{G(\wh{z}_{t+1})}
 \\
 &\le \DD{z{}}{z_{t}} - \frac{\sigma}{2}\twonorm{z_{t+1} - \hz_{t+1}}^2 - 
\frac{\sigma}{2}\twonorm{\hz_{t+1}- z_t}^2 
 + \eta L\twonorm{z_t - \hz_{t+1}}\twonorm{z_{t+1} - \hz_{t+1}}
 \\
 &\qquad\qquad + \eta \iprod{z{} - \hz_{t+1}}{G(\wh{z}_{t+1})}
 \\
 &\le \DD{z{}}{z_{t}} - \frac{\sigma - \eta L}{2}\pa{\twonorm{z_{t+1} - \hz_{t+1}}^2 
+ \twonorm{\hz_{t+1}- z_t}^2}
 \\
 &\qquad\qquad + \eta \iprod{z{} - \hz_{t+1}}{G(\wh{z}_{t+1})}
 \\
 &\le \DD{z{}}{z_{t}} - \frac{\sigma - \eta L}{4}\twonorm{z_{t+1} - z_t}^2 + \eta \iprod{z{} - 
\hz_{t+1}}{G(\wh{z}_{t+1})},
\end{align*}
where we also used the elementary inequalities $2ab\le a^2 + b^2$ and $(a+b)^2 \le 2a^2 + 2b^2$ in 
the last two steps, respectively.
\qed

\subsection{The proof of Lemma~\ref{lem:nubound}}\label{app:nubound}
To enhance readability of the proof, we will omit explicit references to $T$ below, and will simply 
use $\pi$, $\rho$ and $\bmu$ to refer to $\pi_T$, $\rho_T$ and $\bmu_T$, respectively.
Defining $\bnu(x) = \sum_a \bmu(x,a)$ for all $x$, we start by noticing that
\[
\iprod{\bmu}{r} - \rho = \sum_{x,a} \pa{\bnu(x) - \stdist(x)} \pi(a|x) r(x,a) \le 
\onenorm{\bnu - \stdist},
\]
so all we are left with is bounding the total variation distance between $\stdist$ and $\bnu$. To 
do this, we start by fixing an arbitrary $k>0$ and observing that
\begin{equation}\label{eq:nubound}
\begin{split}
\onenorm{\pa{\bnu - \stdist}P_\pi^k} 
 &\le C e^{-k/\tau}\onenorm{\bnu - \stdist}
 \\
 &\le C e^{-k/\tau}\pa{\onenorm{\bnu - \bnu P_\pi^k} + \onenorm{\bnu P_\pi^k - \stdist}},
\end{split}
\end{equation}
where we used Assumption~\ref{ass:mixing} in the first step and the triangle inequality in the 
second one.
Regarding the first term in the parentheses, we repeatedly use the triangle inequality to obtain
\begin{align*}
 \onenorm{\bnu - \bnu P_\pi^k} 
&\le \onenorm{\bnu - \bnu P_\pi}  + \onenorm{\bnu P_\pi - \bnu P_\pi^2} + \dots  + \onenorm{\bnu 
P_\pi^{k-1} - \bnu P_\pi^k}
\\
&= \onenorm{\bnu - \bnu P_\pi}  + \onenorm{\pa{\bnu - \bnu P_\pi} P_\pi} + \dots  + 
\onenorm{\pa{\bnu - \bnu P_\pi} P_\pi^{k-1}}
\\
&\le \onenorm{\bnu - \bnu P_\pi}  + C e^{-1/\tau}\onenorm{\bnu - \bnu P_\pi} + \dots  + 
C e^{-(k-1)/\tau}\onenorm{\bnu - \bnu P_\pi}
\\
&\le C \onenorm{\bnu - \bnu P_\pi}  \sum_{i=0}^{k-1} e^{-i/\tau} \le \frac{C}{1 - e^{-1/\tau}} 
\onenorm{\bnu - \bnu P_\pi}.
\end{align*}
Plugging this bound into Equation~\ref{eq:nubound} and observing that $\bnu P_\pi^k - \stdist = 
\pa{\bnu - \stdist}P_\pi^k$ due to stationarity of $\stdist$, we get
\[
 \onenorm{\pa{\bnu - \stdist}P_\pi^k} \le C e^{-k/\tau}\pa{\frac{C}{1 - e^{-1/\tau}} \onenorm{\bnu 
- 
\bnu P_\pi^k} + \onenorm{\pa{\bnu - \stdist}P_\pi^k}}.
\]
Reordering gives
\[
 \onenorm{\pa{\bnu - \stdist}P_\pi^k} \le \frac{C e^{-k/\tau}}{1 - C e^{-k/\tau}} \cdot \frac{C}{1 
- 
e^{-1/\tau}} \onenorm{\bnu - \bnu P_\pi}.
\]
Thus, using the triangle inequality again yields
\begin{align*}
 \onenorm{\bnu - \stdist} &\le \onenorm{\bnu - \bnu P_\pi^k} + \onenorm{\bnu P_\pi^k - \stdist}
 \\
 &\le \pa{1 + \frac{C e^{-k/\tau}}{1 - C e^{-k/\tau}}}\frac{C}{1 - e^{-1/\tau}} \onenorm{\bnu - 
\bnu 
P_\pi}.
\end{align*}
Now, choosing any $k \ge \tau \log (2C)$ and using the elementary inequality $1/(1-e^{-1/\tau}) 
\le \tau+1$ concludes the proof.
\qed

\subsection{The proof of Lemma~\ref{lem:smoothness_relax}}\label{app:smoothness_relax}
We start by noticing that the 
dual norm of $\norm{z}^2 = \norm{u}_2^2 + \norm{y}_1^2$ evaluated at $x=(w,q)$ is $\norm{x}_*^2 = 
\norm{w}_2^2 + \infnorm{q}^2$. Recalling 
that the smoothness of $\wt{\LL}$ with respect to $\norm{\cdot}$ is equivalent to the Lipschitzness 
of 
$G$ with respect to $\norm{\cdot}_*$, we will prove that $\norm{G(z) - G(z')}_*^2 \le 4K^2\norm{z 
- z'}^2$. 
Using the definition of $G(z)$, we have for any $z=(u,y)$ and $z'=(u',y')$ that
\begin{align*}
 \norm{G(z) - G(z')}_*^2 = \twonorm{F\transpose Q\transpose W\pa{y-y'}}^2 + \infnorm{W\transpose QF(u - 
u')}^2
%  \le \onenorm{Q\transpose W\transpose\pa{y-y'}}^2 + \twonorm{u - u'}^2,
\end{align*}

Let's first see that the sum of any row $j$ of $W\transpose F$ is bounded by $K$:
\[
\sum_i(W\transpose F)_{j,i}=\sum_i(\sum_x W\transpose_{j,x} F_{x,i})= \sum_x W\transpose_{j,x} \sum_i F_{x,i} \leq \sum_x W\transpose_{j,x} K= K 
\]
The same can be easily proven for the matrix $W\transpose PF$.
Now, the first term can be bounded as
\begin{align*}
\twonorm{F\transpose Q\transpose W\pa{y-y'}} &\le \onenorm{F\transpose Q\transpose W\pa{y-y'}} \\
&\leq\onenorm{F\transpose W\pa{y-y'}} + 
\onenorm{F\transpose P \transpose W \pa{y-y'}} \\
& \le 2K\onenorm{y-y'},
\end{align*}
To bound the last term, we observe that
\begin{align*}
\norm{W\transpose QF(u-u')}_\infty^2
& = \max_{j}\left|\sum_{i}\pa{W\transpose F-W\transpose PF}_{j,i} \pa{u_i - u'_i} \right|^2 \\
&\leq \max_{j}\left|\sum_{i} \left(\left|\pa{W\transpose F}_{j,i}\right| +\left|\pa{W\transpose PF}_{j,i}\right|\right) \pa{u_i - u'_i} \right|^2 \\
&\leq \max_{j}\left|\sum_{i} \left(\norm{\pa{W\transpose F}_{j,\cdot}}_1 +\norm{\pa{W\transpose PF}_{j,\cdot}}_1\right) \norm{u - u'}_\infty \right|^2 \\
&\leq 4K^2 \norm{u - u'}_\infty^2 \leq 4K^2 \twonorm{u-u'}^2.
\end{align*}
This concludes the proof.
\qed

\subsection{The proof of Lemma~\ref{lem:yWQbound}}\label{app:yWQbound}
The statement is obvious when $Q\transpose W\transpose\overline{y}_T = 0$, so we will assume 
that the contrary holds below. Let us define
 \[
  w = \taumix \cdot \argmax_{v:\infnorm{v} = 1} \iprod{Q\transpose W\transpose\overline{y}_T}{v},
 \]
 noting that $\iprod{Q\transpose W\transpose \overline{y}_T}{w} = \taumix \onenorm{Q\transpose 
W\transpose \overline{y}_T} > 0$. 
By using this fact and Assumption~\ref{ass:coh}, we crucially observe that there exists a 
$\wt{u}$ such that 
$\iprod{Q\transpose W\transpose \overline{y}_T}{w} = \iprod{Q\transpose W\transpose 
\overline{y}_T}{F\wt{u}}$
and $\infnorm{\wt{u}} \le \taumix U$.  
 This implies that we can apply Corollary~\ref{cor:gapbound} with $z = ( 
F\overline{u}_T  - F\wt{u},  W\transpose \overline{y}_T)$ to obtain the bound
 \[
  \iprod{ Q\transpose W\transpose \overline{y}_T}{w}  = \iprod{ Q\transpose W\transpose 
\overline{y}_T}{ F\overline{u}_T } + \iprod{W\transpose 
\overline{y}_T}{r} - 
\iprod{ Q\transpose W\transpose \overline{y}_T}{ F\pa{\overline{u}_T  - \wt{u}}} - 
\iprod{W\transpose \overline{y}_T}{r}
\le 
\frac{\DD{z}{z_0}}{\eta T}.
 \]
Plugging in the definition of $w$ and the Bregman divergence $D_\Phi$,
we obtain
\begin{align*}
 \onenorm{ Q\transpose W\transpose \overline{y}_T} \le \frac{\frac12 \twonorm{\wt{u} - 
\overline{u}_T}^2  + 
\log 
M}{\eta \taumix T}.
\end{align*}
Due to Assumption~\ref{ass:real} and our assumption on $F$ stated before Theorem~\ref{thm:main}, we 
can choose an optimal solution $u^*$ satisfying $Fu^* = v^*$ and $\infnorm{u^*} \le \taumix U$ and 
write
\begin{align*}
\twonorm{\wt{u}-\overline{u}_T}^2 &\leq 2\twonorm{\wt{u}-u^*}^2+2\twonorm{\overline{u}_T-u^*}^2\leq 
4\twonorm{\wt{u}}^2 + 4\twonorm{u^*}^2+ 4\DD{z^*}{\overline{z}_T} \\
&\leq 4N\infnorm{\wt{u}}^2 + 4N\infnorm{u^*}^2+ 4\DD{z^*}{z_0} \\
& \leq 10\taumix^2 U^2 N + 4 \log M,
\end{align*}
where in the second line we have used Corollary~\ref{cor:normbound} that implies 
$\DD{z^*}{\overline{z}_T}\leq \DD{z^*}{z_0}$.
\qed

\end{document}